\documentclass[11pt,twoside]{article}
\usepackage[colorlinks=true]{hyperref}
\usepackage{pdfsync}
\usepackage [english]{babel}
\usepackage{amsfonts}
\usepackage{amsthm}
\usepackage{amsmath}
\usepackage[hmargin=3cm, vmargin=3cm]{geometry}
\usepackage{amssymb}
\usepackage [T1]{fontenc}
\newtheorem{theorem}{Theorem}[section]
\newtheorem{definition}{Definition}[section]

\newtheorem{remark}{Remark}[section]

\newtheorem{lemma}{Lemma}[section]
\newtheorem{proposition}{Proposition}[section]
\newtheorem{corollary}{Corollary}[section]

\newcommand{\cqfd}
{%
\mbox{}%
\nolinebreak%
\hfill%
\rule{2mm}{2mm}%
\medbreak%
\par%
}
\makeatletter
\newcommand{\cvfstar}{%
  {\mathop{\rightharpoonup}\limits^{\vbox to -\ex@{\kern-\tw@\ex@
        \hbox{\scriptsize *}\vss}}}}
\makeatother
\newcommand{\R}{\mathbb{R}}

\title{\textbf{Insensitizing control for linear and semi-linear heat equations with partially unknown domain}}
\date{}
\author{Pierre Lissy \footnote{Ceremade, Universit\'e Paris-Dauphine \& CNRS UMR 7534, PSL, 75016 Paris, France (\texttt{lissy@ceremade.jussieu.fr}).} \and Yannick Privat\footnote{CNRS, Universit\'e Pierre et Marie Curie (Univ. Paris 6), UMR 7598, Laboratoire Jacques-Louis Lions, F-75005, Paris, France (\texttt{yannick.privat@upmc.fr}).} \and Yacouba Simpor\'e\footnote{Laboratoire LAMI, Université Ouaga 1 Professeur Joseph Ki-Zerbo,
01 BP 7021 Ouaga 01, Burkina Faso (\texttt{simplesaint@gmail.com})}~\footnote{The first author is supported by the ANR project IFSMACS- ANR-15-CE40-0010. The second author is supported by the Project  ``Analysis and simulation of optimal shapes - application
to life science'' of the Paris City Hall.}}
\begin{document}
\maketitle

\begin{abstract}
We consider a semi-linear heat equation with Dirichlet boundary conditions and globally Lipschitz nonlinearity, posed on a bounded domain of $\mathbb R^N$ ($N\in \mathbb N^*$), assumed to be an unknown perturbation of a reference domain. We are interested in an insensitizing control problem, which consists in finding a distributed control such that some functional of the state is insensitive at the first order to the perturbations of the domain. Our first result consists of an approximate insensitization property on the semi-linear heat equation. It rests upon a linearization procedure together with the use of an appropriate fixed point theorem. For the linear case, an appropriate duality theory is developed, so that the problem can be seen as a consequence of well-known unique continuation theorems. Our second result is specific to the linear case. We show a property of exact insensitization for some families of deformation given by one or two parameters. Due to the nonlinearity of the intrinsic control problem, no duality theory is available, so that our proof relies on a geometrical approach and direct computations. 
\end{abstract}
\bigskip \textbf{Keywords}: shape deformation,
insensitizing control, linear and semi-linear heat equation.\\

\noindent \bigskip \textbf{MSC Classification}: 35K05, 35K55, 49K20, 93B05.

\section{Introduction}
\subsection{State of the art}

\indent This article is devoted to proving some results concerning the  insensitizing control for the norm of the linear and semi-linear heat equation when the domain is partially unknown (in the sense that it is a small perturbation of a reference domain).
 The problem of insensitizing control was originally addressed by J-L.Lions in
\cite{16}, leading to numerous papers on this topic.

Concerning the semi-linear heat equation, the first result was obtained in \cite{2} for a distributed control, where the authors introduced and studied the
notion of approximate insensitizing controls, the partially unknown data being the initial condition and the boundary condition, the sentinel being the square of the $L^2-$norm of the solution on some subset of $\Omega$, called observation domain. The method used in this article (that has widely inspired the present article) is to study first a linear heat equation with potential and then apply the Schauder's fixed point
theorem to extend the conclusions to the semi-linear heat equation. This result was improved in \cite{1}, where the existence of exact insensitizing controls for perturbations of the null initial datum is proved by means of Carleman estimates (see \cite{FI}) for a linear forward-backward coupled heat system and the use of a fixed-point argument. Let us also mention that in \cite{1}, it is proved that when the initial datum is not null, one cannot always expect to find an insensitizing control (for a precise study of the class of initial data that can be insensitized for the linear heat equation, see \cite{5}).  Later on, the insensitizing control problem for the heat equation with nonlinear boundary conditions was studied notably in \cite{3}. In \cite{Gu}, an exact insensitizing control result for  a linear heat equation with potential is proved  when the sentinel  is the square of the $L^2$-norm of the gradient of the solution on some observation subset. In \cite{20}, the authors studied the insensitizing control problem with constraints on the control for a nonlinear heat equation by means of the Kakutani's fixed point theorem combined with an adapted Carleman inequality. This result was extended in \cite{SO} for more general cost functionals. Some quasi-linear parabolic problems have also  been studied in \cite{Liu}.

Let us mention that in all the articles mentioned above, a crucial hypothesis is that the observation domain intersects the control domain. Removing this hypothesis leads to many difficulties (notably because Carleman estimates cannot be used anymore),  this case being notably studied in \cite{Mi} and \cite{KL}.

To conclude, let us mention that other linear or nonlinear parabolic systems coming from fluid mechanics have also been intensively studied, see for instance \cite{Gu2}, \cite{Gueye}, \cite{CG}, \cite{CGG} or \cite{Car}.

The common point of all the previous articles is that the partially known datum considered is always the initial condition, with sometimes the addition of the boundary condition. Up to our knowledge, the question of insensitizing controls for a deformation of the domain has never been studied. Let us mention that a close problem was studied in \cite{BD}, but the goal of the authors was different since they intended to estimate the shape of an unknown part of a domain for a diffusion problem. In the framework of control theory, we also mention \cite{PS,MPS} where genericity of controllability properties with respect to domain variations are investigated.

The general problem of trying to insensitizing an observation done on a domain that is partially unknown is meaningful from the applicative point of view. One can for example think of an oil drilling: we observe the drilling  on a known domain, the initial shape of the oil field is known, but the extraction may perturbe the shape of this field. The goal may be then to optimize the observation  by acting on some other place of the field. 

The rest of this paper is organized as follows: in Section \ref{state} we present the problem and our main results. In Section \ref{CO}, we use standard arguments to reduce the problem to a control problem on a forward-backward system of semi-linear heat equations. Section \ref{S1} is devoted to the proof of Theorem \ref{th1}, whereas Section \ref{S2} is devoted to the proof of Theorem \ref{th2}.  To finish, we give some perspectives in the concluding Section \ref{CON}.

 \subsection{Statement of the problem and main result}\label{state}
Let $\Omega_{0}$ be a connected and bounded subset of $\mathbb{R}^{N},\text{ } N\in\mathbb N^*,$ assumed to be of class $C^2$. Let $T>0$ and let $\omega$
and $\Theta$ be two nonempty open subsets of $\Omega_{0}$, assumed to be compactly included in $\Omega_0$. We set
 $Q_{0}=\left(0,T\right)\times\Omega_{0}$
 and $\Sigma_0=\left(0,T\right)\times\partial\Omega_{0}$.

Since we are interested in dealing with perturbations of $\Omega_0$ preserving some topological properties such as its connectedness, boundedness and regularity, we will adopt the classical point of view in shape optimization used to define the derivative in the sense of  Hadamard (see e.g. \cite{DZ,Henrot-Pierre}). This means that perturbations of $\Omega_0$ will be defined with the help of well-chosen diffeomorphisms.

In this view, let us introduce for any integer $j\geq 1$ the admissible class of perturbations fields
$$
\mathcal{V}^{j,\infty}=\{\mathbf{V} \in W^{j,\infty}(\R^N,\R^N) \mid  \Vert\mathbf{V}\Vert_{j,\infty}\leqslant 1\}.
$$
It is notable that, for each element $\mathbf{V}$ of $\mathcal{V}^{3,\infty}$ and each $\tau\in [0,1)$, the mapping $\mathbf{T_\tau}:=\operatorname{Id}+\tau\mathbf{V}$ defines a diffeomorphism in $\R^d$, i.e. the mapping $\mathbf{T_\tau}$ is invertible and $\mathbf{T_\tau}^{-1}\in W^{3,\infty}(\R^N,\R^N)$.
Furthermore, as a consequence of the construction of $\mathbf{T_\tau}$ as a ``perturbation of the identity'', the set $\mathbf{T_\tau}(\Omega_0)$ is a connected, bounded domain whose boundary is of class $C^2$. 

In the sequel, we will consider a family of domains $\{\Omega_\tau\}_{\tau\in [0,1)}$ of $\Omega_0$ defined, for some $\mathbf{V}\in \mathcal{V}^{3,\infty}$, by 
$$
\Omega_\tau= (\operatorname{Id}+\tau\mathbf{V})(\Omega_0).
$$
As a consequence, each domain $\Omega_\tau$ inherits the aforementioned properties, moreover the sets $\omega$ and $\Theta$ are compactly included in $\Omega_\tau$ provided that the parameter $\tau$ is chosen small enough, which is assumed from now on.

 Let us set
 $Q_{\tau}=\left(0,T\right)\times\Omega_{\tau}$
 and $\Sigma=\left(0,T\right)\times\partial\Omega_{\tau}$. Let $\chi_{\omega}\text{ and }\chi_{\Theta}$ respectively be the characteristic function of $\omega\text{ and } \Theta$.\\
This article is concerned with the family of systems
\begin{equation}\label{nlh}
\left\{ 
\begin{array}{rcll}
\dfrac{\partial y}{\partial t}-\Delta y +f(y)&=&\xi +h\chi _{\omega
}& \text{ in }Q_{\tau} ,\\ 
y&=& 0 & \text{ on }\Sigma_{\tau},\\ 
y\left( 0,\cdot\right) &=& 0 & \text{ in }\Omega_{\tau},
\end{array}%
\right.
\end{equation}
where $f\in C^1(\mathbb R)$ is assumed to be globally  Lipschitz and  $\xi\in L^{2}( \mathbb R^N) $. The control term $h$ belongs to $L^{2}((0,T)\times\omega).$
The data of the state equation \eqref{nlh} are incomplete in the sense that both $\mathbf{V}$ and $\tau$ are partially unknown.

\begin{remark}
Let us comment on the initial condition. For the sake of simplicity, we chose to consider an identically null initial condition. However, all the results of this article are easily generalizable to any initial  condition $y(0,\cdot)=y^0$ where $y^{0}\in H^{1}_0(\Omega_0)$ is such that $\operatorname{dist}(\operatorname{supp}(y^0),\partial\Omega_0)>0$. Indeed, the $H^1$ regularity of $y^0$ is needed to use the insensitizing conditions (see Remark \ref{rk:smooth}) whereas the distance condition on the support of $y^0$ guarantees that the initial boundary condition will not be sensitive to the variations of $\Omega_0$. 
\end{remark}

Let us now provide a precise definition of insensitizing a functional with respect to a deformation of the domain.

\begin{definition}\label{def:eiai}
Let $\Phi:L^2(0,T,L^2(\R^N))\ni y \mapsto \Phi(y)\in \R$ be a differentiable functional and let $\Omega$ be a connected bounded domain having a $C^2$ boundary. Let us introduce the shape functional $\mathcal{J}$ defined by $\mathcal{J}(\Omega)=\Phi(y_{\Omega})$ where $y_{\Omega}$ denotes the extension to $0$ on $\mathbb R_+^*\times \mathbb{R}^N$ of the unique solution of $\eqref{nlh}$, where $\Omega_\tau$ has been replaced by $\Omega$.

For some $\mathbf{V}\in \mathcal{V}^{3,\infty}$, we introduce as above the domains $\Omega_\tau= (\operatorname{Id}+\tau\mathbf{V})(\Omega_0)$ (where the parameter $\tau$ lies  in $[0,\tau_0)$ for some small enough $\tau_0<1$).
Let $\xi \in L^2(\R^N)$ be given. We say that the control $h$ insensitizes $\Phi$ (at the first order) whenever
\begin{equation}\label{ei}
	\text{for all } \mathbf{V}\in \mathcal{V}^{3,\infty},\quad \text{there holds } \left.\dfrac{d }{d \tau}\left(\mathcal{J}\left(\Omega_\tau\right)\right)\right|_{\tau=0}=0.
\end{equation}
Let $\mathcal{E}$ be a linear subspace of $\mathcal{V}^{3,\infty}$. One says that the control $h$ insensitizes $\Phi$ (at the first order) for the family $\mathcal{E}$ whenever
\begin{equation}\label{eiFinite}
	\text{for all } \mathbf{V}\in \mathcal{E},\quad \text{there holds } \left.\dfrac{d }{d \tau}\left(\mathcal{J}\left(\Omega_\tau\right)\right)\right|_{\tau=0}=0.
\end{equation}

Given $\varepsilon>0,$ the control $h$ is said to $\varepsilon$-insensitize $\Phi$ whenever
\begin{equation}\label{ai}
\text{for all }  \mathbf{V}\in \mathcal{V}^{3,\infty},\quad \text{there holds } \left|\left.\dfrac{d }{d \tau}\left(\mathcal{J}\left(\Omega_\tau\right)\right)\right|_{\tau=0}\right |\leq \varepsilon.
\end{equation}
\end{definition}

Notice that this definition uses a particular notion of derivative, well adapted when dealing with shape variations. In this setting, we consider variations of a domain that are parametrized by families of diffeomorphisms, as highlighted previously. 

In what follows, we will concentrate on a particular choice of shape functional $\mathcal{J}$ that appears natural in the framework of control: we aim at insensitizing the $L^2$-norm of the solution $y_\Omega$ of \eqref{nlh} with respect to the domain, which leads to consider the functional  
\begin{equation}\label{defphi}
\mathcal{J}(\Omega)=\Phi(y_\Omega)=\dfrac{1}{2}\int\limits_{0}^{T}\!\!\! \int_{\Theta}y_\Omega(t,x)^{2}\, dxdt.
\end{equation}
In other words, $\mathcal{J}$ stands for the square of the $L^{2}$ norm of the observation variable $\chi_\Theta y_\Omega$.

Let us now describe the main results of this article.

\begin{theorem}\label{th1}Assume that $\omega\cap\Theta\not = \emptyset$. Then, for every $\varepsilon>0$, there exists a control $h\in L^2((0,T)\times\omega)$ which $\varepsilon$-insensitizes $\Phi$.
\end{theorem}

\begin{theorem}\label{th2}
Let $M\in \{1,2\}$. Assume that $\omega\cap \Theta\neq \emptyset$ and let $\{\mathbf{V}_i\}_{1\leq i\leq M}$ be a family of linearly independent elements of $\mathcal{V}^{3,\infty}$. Then, there exists a control $h$ which insensitizes the functional $\Phi$ for the family $\mathcal{E}=\operatorname{span}(\{\mathbf{V}_i\}_{1\leq i\leq M})$.
\end{theorem}
Several geometrical examples of framework and families $\mathcal{E}$ for which Theorem \ref{th2} applies are provided in Section \ref{sec:settingExactInsens}. 
\begin{remark}[Comments on the regularity of $\Omega$]\label{rk:smooth}
Notice that the notion of shape derivative does not impose to deal with regular shapes. Indeed, considering a bounded connected domain $\Omega_0$ with a Lipschitz boundary is enough to define the quantities \eqref{ei} and \eqref{ai} involved in Definition \ref{def:eiai}. In particular, in that setting, the mapping $\mathcal{V}^{1,\infty}\ni \mathbf{V}\mapsto \mathcal{J}((\operatorname{Id}+\mathbf{V})(\Omega_0))\in \R$ is differentiable and there holds in particular
$$
\langle d\mathcal{J}(\Omega),\mathbf{V}\rangle =\lim_{\tau \searrow 0}\frac{\mathcal{J}(\Omega_\tau)-\mathcal{J}(\Omega_0)}{\tau}.
$$
Nevertheless, although well-defined, the differential of $\mathcal{J}$ at $\Omega$ cannot be recast without additional regularity assumptions in a simple form and is difficult to handle. This is why we chose to deal with domains having a $C^2$ boundary (and then to consider perturbations in $\mathcal{V}^{3,\infty}$).
\end{remark}
\subsection{Reduction of the insensitizing control to a control problem on a coupled system}
\label{CO}
Let us consider a domain $\Omega_0$ enjoying the same properties as in Section \ref{state}. This section is devoted to deriving insensitizing conditions on the domain $\Omega_0$, in other words to recast the conditions \eqref{ai} and \eqref{ei} in a simpler way. We will emphasize that the problem of exact (resp. approximate) insensitizing control can be reduced to a non-standard null (resp. approximate) controllability problem on a backward-forward coupled system of semi-linear heat equations. We first claim that the mapping
$$
\mathcal{V}^{3,\infty}\ni \mathbf{V}\mapsto y_{(\operatorname{Id}+\mathbf{V})(\Omega_0)}\in L^2(0,T,L^2(\R^N))
$$
is differentiable at $\mathbf{V}=0$, and therefore, so is the mapping $\R \ni \tau \mapsto y_{\Omega_\tau} \in L^2(0,T,L^2(\R^N))$ at $\tau=0$. The reasoning to obtain this result is standard and rests upon the implicit function theorem (see e.g. \cite[Theorem 5.3.2]{Henrot-Pierre}). We denote by $\dot y_{\Omega_0}$ the differential of this mapping (also called the Eulerian derivative of $y_{\Omega_0}$) in a given direction $\mathbf{V}$.

It is well-known (see e.g. \cite[Theorem 5.3.1]{Henrot-Pierre} or also \cite{Soko}) that the function $\dot y_{\Omega_0}$ solves the partial differential equation
\begin{equation}\label{htau}
\left\{ 
\begin{array}{rclc}
\dfrac{\partial \dot y_{\Omega_0}}{\partial t}-\Delta \dot y_{\Omega_0} +f'(y_{\Omega_0})\dot y_{\Omega_0} &=&0 & \text{ in } Q_{0} ,\\ 
\dot y_{\Omega_0} &=& -\partial_{n}y_{\Omega_0} (\mathbf{V}\cdot \mathbf{n}) & \text{ on }\Sigma_{0},\\ 
\dot y_{\Omega_0}\left( 0,\cdot \right) &=& 0 & \text{ in }\Omega_{0},
\end{array}%
\right.
\end{equation}
where $\mathbf{n}$ stands for the unit outward normal of $\partial \Omega_0$.
According to \eqref{ei}-\eqref{ai}-\eqref{defphi}, the exact insensitization control problem comes to 
 \begin{equation}\label{newe}
	\forall \mathbf{V} \in \mathcal{V}^{3,\infty},\quad \int_{0}^{T}\!\!\! \int_{\Theta}y_{\Omega_0}(t,x)\dot y_{\Omega_0}\,  dx dt=0,
\end{equation}
whereas the $\varepsilon$-approximate insensitization control problem is equivalent to 
 \begin{equation}\label{newa}
	\forall \mathbf{V} \in \mathcal{V}^{3,\infty},\quad \left|\int_0^T\int_{\Theta}y_{\Omega_0}(t,x)\dot y_{\Omega_0}\,  dx dt \right|\leqslant \varepsilon.
\end{equation}
We are going to provide a more workable characterization of (exact and approximate) insensitizing conditions in terms of the solutions of a forward-backward coupled system.
\begin{proposition}\label{propcou}
The exact insensitizing control problem \eqref{ei} is equivalent to the following one: for any $\xi\in L^2(Q_0)$, find $h\in L^2(Q_0)$ such that the solution $(y_{\Omega_0},q_{\Omega_0})$ to the following forward-backward coupled system
	\begin{equation}
\left\{ 
\begin{array}{rcll}
\dfrac{\partial y}{\partial t}-\Delta y +f(y)&=&\xi +h\chi _{\omega} & \text{ in }Q_{0} ,\\ 
y&=& 0 & \text{ on }\Sigma_{0},\\ 
y\left( 0,\cdot \right) &=& 0 & \text{ in }\Omega_{0},
\end{array}%
\right.\label{cou1}
\end{equation}
\begin{equation}
\left\{ 
\begin{array}{rcll}
-\dfrac{\partial q}{\partial t}-\Delta q +f'(y)q &=& y\chi_{\Theta}& \text{ in }Q_{0} ,\\ 
q&=& 0 & \text{ on }\Sigma_{0},\\ 
q\left( T,\cdot\right) &=& 0 & \text{ in }\Omega_{0},
\end{array}%
\right.\label{cou2}
\end{equation}
satisfies
\begin{equation}\label{efin}
	\int_{0}^{T}\partial_{n} y_{\Omega_0}\partial_{n} q_{\Omega_0}\, dt=0, \quad \mbox{a.e. in } \partial \Omega_{0}.
\end{equation}

The $\varepsilon-$approximate insensitizing problem \eqref{ei} is equivalent to the following problem: for any $\xi\in L^2(Q_0)$, find $h\in L^2(Q_0)$ such that the solution $(y_{\Omega_0},q_{\Omega_0})$ to \eqref{cou1}-\eqref{cou2} satisfies 

\begin{equation}\label{afin}
	\int_{\partial\Omega_0}\left |\int_{0}^{T}\partial_{n} y_{\Omega_0}\partial_{n} q_{\Omega_0}\, dt  \right |d\sigma\leqslant \varepsilon.
\end{equation}
\end{proposition}

Notice that the resulting control problems \eqref{efin} and \eqref{afin} on the solutions of \eqref{cou1}-\eqref{cou2} is quite unusual since it writes as a \emph{ bilinear problem} with respect to the states $y_{\Omega_0}$ and $q_{\Omega_0}$, which makes its study more difficult than for standard problems of controllability. Notably, even if we consider linear versions of the system \eqref{cou1}-\eqref{cou2}, the control problems \eqref{efin} and \eqref{afin} are bilinear, so that the standard duality theory (see \cite{LHUM} or \cite{TW} for instance) cannot be applied.

\vspace*{0.3cm}

\noindent \textbf{Proof of Proposition \ref{propcou}.}

Let us first remark that since $(\xi,h\chi_\omega)\in L^2(Q_0)\times L^2(Q_0)$ and $y(0,\cdot)=0$, the solution $y_{\Omega_0}$ of \eqref{cou1} satisfies 
$$
y_{\Omega_0}\in L^2(0,T,H^1_0(\Omega_0)\cap H^2(\Omega_0))\cap H^1(0,T,L^2(\Omega_0)).
$$
As a consequence, the function $\partial_n y_{\Omega_0}$ is well-defined on $\Sigma_0$ and we claim moreover that $\partial_n y_{\Omega_0}$ belongs to $ L^2(0,T,L^2(\partial\Omega_0))$. The same argument enables us to show that the solution $q_{\Omega_0}$ of \eqref{cou2} satisfies 
 $$
 q_{\Omega_0}\in L^2(0,T,H^1_0(\Omega_0)\cap H^2(\Omega_0))\cap H^1(0,T,L^2(\Omega_0)),
 $$
 so that $\partial_n q_{\Omega_0}$ also makes sense on $\Sigma_0$. Moreover, one has $\partial_n q_{\Omega_0}\in L^2(0,T,L^2(\partial\Omega_0))$. 
 
 We infer notably that the mapping
$$
\partial\Omega_0\ni x\mapsto \int_{0}^{T}\partial_{n} y_{\Omega_0}(t,x)\partial_{n} q_{\Omega_0}(t,x)\,  dt 
$$
belongs to $L^1(\partial\Omega_0)$.

Multiplying the first equation of \eqref{cou2} by $\dot y_{\Omega_0}$ and integrating by parts yields 
\begin{eqnarray*}
\lefteqn{-\int_0^T\!\!\! \int_{\Omega_0}\frac{\partial q_{\Omega_0}}{\partial t}\dot y_{\Omega_0}\, dt dx+\int_0^T\!\!\! \int_{\Omega_0} (\nabla q_{\Omega_0}\cdot \nabla \dot y_{\Omega_0}+f'(y_{\Omega_0})q_{\Omega_0}\dot y_{\Omega_0})\, dxdt } \phantom{aaaaaaaaaaaaaaaaa}\\
& & \displaystyle  +\int_0^T\!\!\!\int_{\partial\Omega_0}\frac{\partial q_{\Omega_0}}{\partial n}\frac{\partial y_{\Omega_0}}{\partial n}(\mathbf{V}\cdot \mathbf{n})\, d\sigma dt= \int_0^T\!\!\! \int_{\Theta} \dot y_{\Omega_0}y_{\Omega_0}\, dx dt. 
\end{eqnarray*}
Similarly, multiplying the first equation of \eqref{htau} by $q_{\Omega_0}$ and integrating by parts yields 
$$
\int_0^T\!\!\! \int_{\Omega_0}\frac{\partial \dot y_{\Omega_0}}{\partial t}q_{\Omega_0}\, dt dx+\int_0^T\!\!\! \int_{\Omega_0} (\nabla q_{\Omega_0}\cdot \nabla \dot y_{\Omega_0}+f'(y_{\Omega_0})q_{\Omega_0}\dot y_{\Omega_0})\, dxdt = 0. 
$$
The combination of the two last equalities leads to
\begin{equation}\label{ci1}
	\int_{0}^{T}\!\!\!\int_{\partial\Omega_{0}}\partial_{n} y_{\Omega_0}\partial_{n} q_{\Omega_0}(\mathbf{V}\cdot \mathbf{n}) \, d\sigma  dt=\int_{0}^{T}\!\!\!\int_{\Theta}y(0)\dot y_{\Omega_0}\, dxdt .
\end{equation}
As a consequence, it follows that \eqref{newe} is equivalent to 
\begin{equation*}
	\int_{0}^{T}\!\!\! \int_{\partial\Omega_{0}}\partial_{n} y_{\Omega_0}\partial_{n} q_{\Omega_0} (\mathbf{V}\cdot \mathbf{n}) \, d\sigma  dt=0,\, \quad \forall \mathbf{V}\in\mathcal{V}^{3,\infty},
\end{equation*}
which rewrites also
\begin{equation*}
	\int_{\partial\Omega_{0}}(\mathbf{V}\cdot \mathbf{n})\left(\int_{0}^{T}\partial_{n} y_{\Omega_0}\partial_{n} q_{\Omega_0}\, dt \right)\, d\sigma=0,\, \quad \forall \mathbf{V}\in\mathcal{V}^{3,\infty}.
\end{equation*}
This equality is equivalent by density and linearity to
\begin{equation*}
	\int_{\partial\Omega_{0}}\beta(x)\left(\int_{0}^{T}\partial_{n} y_{\Omega_0}\partial_{n} q_{\Omega_0} \, dt \right)\, d\sigma=0,\, \quad \forall \beta\in C^2(\partial\Omega_0).
\end{equation*}
Since the mapping $\partial\Omega_0\ni x \mapsto \int_{0}^{T}\partial_{n} y_{\Omega_0}(t,x)\partial_{n} q_{\Omega_0}(t,x)\, dt$ belongs to $L^1(\partial\Omega_0)$, we conclude by applying the fundamental lemma of calculus of variations that the previous equality is equivalent to \eqref{efin}.

On the other hand, using \eqref{ci1}, we know that \eqref{newa} is equivalent to 

\begin{equation*}
	\left| \int_{0}^{T}\!\!\! \int_{\partial\Omega_{0}}\partial_{n} y_{\Omega_0}\partial_{n} q_{\Omega_0} (\mathbf{V}\cdot \mathbf{n}) \, d\sigma  dt \right |\leqslant \varepsilon, \quad  \forall \mathbf{V}\in\mathcal{V}^{3,\infty},
\end{equation*}
i.e.
\begin{equation*}
	\left | \int_{\partial\Omega_{0}}(\mathbf{V}\cdot \mathbf{n})\left(\int_{0}^{T}\partial_{n} y_{\Omega_0}\partial_{n} q_{\Omega_0}\, dt \right)\, d\sigma\right |\leqslant \varepsilon,  \quad  \forall \mathbf{V}\in\mathcal{V}^{3,\infty},
\end{equation*}
This inequality  is equivalent by density and linearity to
\begin{equation*}
\left | \int_{\partial\Omega_{0}}\beta(x)\left(\int_{0}^{T}\partial_{n} y_{\Omega_0}\partial_{n} q_{\Omega_0} \, dt \right)d\sigma\right |\leqslant \varepsilon ||\beta||_\infty, \quad  \forall \beta\in C^2(\partial\Omega_0),
\end{equation*}
which is also equivalent by density to
\begin{equation*}
\left | \int_{\partial\Omega_{0}}\beta(x)\left(\int_{0}^{T}\partial_{n} y_{\Omega_0} \partial_{n} q_{\Omega_0}\, dt \right)d\sigma\right |\leqslant \varepsilon ||\beta||_\infty, \quad \forall \beta\in L^\infty(\partial\Omega_0),
\end{equation*}
By duality, this exactly means \eqref{afin} and the proof is complete.
\cqfd

\section{Approximated null controllability and approximated sentinel}\label{S1}
\subsection{An auxiliary linear problem}

In a first time, let us investigate  in details the following forward-backward coupled system of linear equations:
\begin{equation}
\left\{ 
\begin{array}{rcll}
\dfrac{\partial u}{\partial t}-\Delta u +a(t,x)u&=&\xi +k\chi _{\omega} & \text{ in }Q_{0} ,\\ 
u&=& 0 & \text{ on }\Sigma_{0},\\ 
u\left( 0,.\right) &=& 0 & \text{ in }\Omega_{0},
\end{array}%
\right.\label{lh1}
\end{equation}
\begin{equation}
\left\{ 
\begin{array}{rcll}
-\dfrac{\partial v}{\partial t}-\Delta v +b(t,x)v &=& u\chi_{\Theta} & \text{ in }Q_{0} ,\\ 
v&=& 0 & \text{ on }\Sigma_{0},\\ 
v\left( T,.\right) &=& 0 & \text{ in }\Omega_{0},
\end{array}%
\right.\label{lh2}
\end{equation}
where $a,b\in L^\infty(Q_0)$ and $\xi\in L^2(Q_0)$.
Before giving the main result of this part, let us provide several useful duality results.
We introduce the operators 
\begin{equation}
\label{deffghi}
\begin{array}{c}
 \mathcal F: L^2(Q_0)\ni \xi \longmapsto (\partial_n u_\xi ,\partial_n v_\xi ) \in L^2(\Sigma_0)\times L^2(\Sigma_0),\\
 \mathcal G: L^2(Q_0)\ni k\longmapsto  (\partial_nu_k,\partial_nv_k) \in L^2(\Sigma_0)\times L^2(\Sigma_0),\\
\end{array}
\end{equation}
where the pair $(u_\xi, v_\xi)$ denotes the unique solution of the coupled system \eqref{lh1}-\eqref{lh2} in the case $k=0$, and the pair $(u_\xi, v_\xi)$ denotes the unique solution of the coupled system \eqref{lh1}-\eqref{lh2} in the case $\xi=0$.

\medskip

Let us compute the adjoint operator of $\mathcal G$.  
For that purpose, let us consider $(\delta_1,\delta_2)\in (L^2(\Sigma_0))^2$ as well as $(\varphi,\psi)$ solving the system 
\begin{equation}
\left\{ 
\begin{array}{rcll}
\dfrac{\partial \varphi}{\partial t}-\Delta \varphi +b(t,x)\varphi&=&0 & \text{ in }Q_{0} ,\\ 
\varphi &=&\delta_2& \text{ on }\Sigma_{0},\\ 
\varphi\left( 0,.\right) &=& 0& \text{ in }\Omega_{0},
\end{array}%
\right.\label{adj1}
\end{equation}
\begin{equation}
\left\{ 
\begin{array}{rcll}
-\dfrac{\partial \psi}{\partial t}-\Delta \psi +a(t,x)\psi &=& \varphi\chi_{\Theta} & \text{ in }Q_{0} ,\\ 
\psi&=&  \delta_1 & \text{ on }\Sigma_{0},\\ 
\psi\left( T,.\right) &=& 0& \text{ in }\Omega_{0},
\end{array}%
\right.\label{adj2}
\end{equation}
Then, there holds
\begin{eqnarray*}
\langle \mathcal G(k),(\delta_1,\delta_2)\rangle _{L^2(\Sigma_0)\times L^2(\Sigma_0)}&= & \int_{\Sigma_0} \partial_nu \delta_1+\int_{\Sigma_0} \partial_nv \delta_2  \\ 
&=& \int_{Q_0} (\Delta u)\psi-\int_{Q_0} u(\Delta \psi)+\int_{Q_0} (\Delta v)\varphi-\int_{Q_0} v(\Delta \varphi)\\ 
&=& \int_{Q_0} (\partial_t u+a(t,x)u-k\chi_{\omega})\psi-\int_{Q_0} u(-\partial_t \psi+a(t,x)\psi-\varphi\chi_\Theta)\\
& & +\int_{Q_0} (-\partial_t v+b(t,x)v-u\chi_\Theta)\varphi-\int_{Q_0} v(\partial_t \varphi+b(t,x)\varphi)\\
&=& -\int_{Q_0}k\chi_{\omega}\psi.
\end{eqnarray*}
We then infer that
\begin{equation}\label{adjG}
\mathcal G^*: L^2(\Sigma_0)\times L^2(\Sigma_0) \ni (\delta_1,\delta_2) \longmapsto - \psi \chi_\omega \in L^2(Q_0),\end{equation}
where $\psi$ is the solution of \eqref{adj2}.

The end of this section is devoted to introducing a constructive approach for building $\varepsilon$-insensitizing controls in the linear case.

\begin{proposition}\label{mainpr}
Assume that $\omega\cap \Theta\neq \emptyset$. Let $\varepsilon>0$, $\xi\in L^2(\Omega)$, and $(\gamma_1,\gamma_2)\in L^2(\Sigma_0)\times L^2(\Sigma_0)$. There exists $k\in L^{2}\left((0,T)\times \omega\right)$ such that the solution $(u_{\Omega_0},v_{\Omega_0})$ of \eqref{lh1}-\eqref{lh2} satisfies 
	\begin{equation}
\left .
\begin{array}{c}\label{linep}
\displaystyle\left\Vert\partial_{n} u_{\Omega_0}-\gamma_1   \right \Vert_{L^2(\Sigma_0)}^2+\left\Vert\partial_{n} v_{\Omega_0}-\gamma_2   \right \Vert_{L^2(\Sigma_0)}^2\leqslant \varepsilon^2.
\end{array}
\right .
\end{equation}
Let us introduce
\begin{equation*}
\mathcal U_{ad}(T,\varepsilon):=\{k\in L^2(Q_0)\text{ s.t. the solution $(u_{\Omega_0},v_{\Omega_0})$ of \eqref{lh1}-\eqref{lh2} satisfies }\eqref{linep}\},
\end{equation*}
as well as the cost functional $J_\varepsilon$ defined on $(L^2(\Sigma_0))^2$ by
\begin{equation*}
J_\varepsilon(\delta_1,\delta_2):=\frac{1}{2}\int_{(0,T)\times\omega} \psi^2+\varepsilon\displaystyle\sqrt{\Vert \delta_1\Vert ^2_{L^2(\Sigma_0)}+\Vert \delta_2\Vert ^2_{L^2(\Sigma_0)}}-\int_{\Sigma_0}\gamma_1\delta_1+\int_{\Sigma_0}\gamma_2\delta_2.
\end{equation*}
\begin{enumerate}
\item There holds
$$\min_{k\in\mathcal U_{ad}(T,\varepsilon)} \frac{1}{2}\int_{(0,T)\times\omega}k^2=-\min_{(\delta_1,\delta_2)\in L^2(\Sigma_0)\times L^2(\Sigma_0)}J_\varepsilon(\delta_1,\delta_2).$$
\item The control $k_{opt}^\varepsilon$ of minimal $L^2$-norm is given by $k_{opt}^\varepsilon=\psi_{opt}^\varepsilon \chi_\omega$, where $\psi_{opt}^\varepsilon$ is the solution of \eqref{adj2} associated to the minimum $(\delta_{1,opt}^\varepsilon ,\delta_{2,opt}^\varepsilon)$ of $J_\varepsilon$.
\item If $a$ and $b$ describe bounded sets of $L^\infty$, if $\gamma_1$ and $\gamma_2$ describe compact sets of $L^2(\Sigma)$, then the control obtained above describes a  bounded (and even compact) set of $L^2(Q_0)$.

\end{enumerate}
\end{proposition}

Notice that, by applying the Cauchy-Schwarz inequality together with Young's inequality, Proposition \ref{mainpr} (with $\gamma_1=\gamma_2=0$) implies that for any $\varepsilon>0$, there exists $h_\varepsilon \in L^{2}\left((0,T)\times \omega\right)$ such that the solution $(u_{\Omega_0},v_{\Omega_0})$ of \eqref{lh1}-\eqref{lh2} with $h=h_\varepsilon$ satisfies 

\begin{equation*}
	\int_{\partial\Omega_0}\left |\int_{0}^{T} \partial_{n} u_{\Omega_0}\partial_{n} v_{\Omega_0} \, dt \right | d\sigma\leqslant \varepsilon,
\end{equation*}
in accordance with the inequality \eqref{afin} in Section \ref{CO}.

\medskip

\noindent \textbf{Proof of Proposition \ref{mainpr}.}
Using \eqref{deffghi}, we observe that \eqref{linep} is equivalent to asking that for for any $\varepsilon>0$, any $\xi\in L^2(\Omega)$,  and any $(\gamma_1,\gamma_2)\in L^2(\Sigma_0)\times L^2(\Sigma_0)$, one has
$$
\Vert \mathcal F(\xi)+\mathcal G(k)-(\gamma_1,\gamma_2)\Vert_{L^2(\Sigma_0)\times L^2(\Sigma_0)} \leqslant \varepsilon.
$$
Hence, the property ``for any $\varepsilon>0$, any $\xi\in L^2(\Omega)$ and any  $(\gamma_1,\gamma_2)\in L^2(\Sigma_0)\times L^2(\Sigma_0)$, there exists $h\in L^{2}\left((0,T)\times \omega\right)$ such that the solution $(u,v)$ of \eqref{lh1}-\eqref{lh2} satisfies \eqref{linep}'' will be true as soon as we are able to prove that
$$\overline{\text{Range }(\mathcal G)}=L^2(\Sigma_0)\times L^2(\Sigma_0).$$
By duality, this is also equivalent to the following unique continuation property:
$$\mathcal G^*(\delta_1,\delta_2)=0\Rightarrow \delta_1=\delta_2=0.$$
Using \eqref{adjG}, this property also rewrites as 

\begin{equation}\label{ucppp}\psi=0 \text{ in } (0,T)\times\omega\Rightarrow \delta_1=\delta_2=0 \text{ in } L^2(\Sigma_0),
\end{equation}
where $\psi$ denotes the solution of \eqref{adj2}. 

Let us prove property \eqref{ucppp}.
Assume that $\psi=0$ on $(0,T)\times\omega$. Then,  using \eqref{adj2}, we infer that $\varphi=0$ on $(0,T)\times(\Theta\cap\omega)$, where $\varphi$ is the solution of \eqref{adj1}. Recall that one has $\Theta\cap\omega\not=\emptyset$ by assumption. Hence, applying \cite[Proposition 2, Page 670]{2} (which is a consequence of the unique continuation results of \cite{SS}), we infer that $\varphi=0$ on $Q_0$ and $\delta_2=0$ on $\Sigma_0$. Hence, $\psi$ satisfies the backward equation

\begin{equation*}
\left\{ 
\begin{array}{rcll}
-\dfrac{\partial \psi}{\partial t}-\Delta \psi +a(t,x)\psi &=&0 & \text{ in }Q_{0} ,\\ 
\psi&=& \delta_1 & \text{ on }\Sigma_{0},\\ 
\psi\left( T,.\right) &=& 0 & \text{ in }\Omega_{0},
\end{array}%
\right.
\end{equation*}
together with $\psi=0 \text{ in } (0,T)\times\omega$, so that we can one more time apply \cite[Proposition 2, Page 670]{2} and deduce that one has necessarily $\delta_1=0$ on $\Sigma_0$.

Points (i) and (ii) are very classical duality results that can be interpreted as consequences of the Fenchel-Rockafellar theory. They may be obtained by applying the so-called Hilbert Uniqueness Method (HUM) and are left to the reader (see for example \cite[Proof of Theorem 3]{2}).

Finally, the last point (iii) is less classical but may be recovered by following the approach in \cite[Proof of Theorem 3]{2}. 
Let us provide hereafter a complete proof for the sake of completeness.

Assume that $a$ and $b$ describe bounded sets of $L^\infty$, whereas $\gamma_1$ and $\gamma_2$ describe compact sets of $L^2(\Sigma)$.
Let us first prove that $(\delta_{1,opt}^\varepsilon,\delta_{2,opt}^\varepsilon)$ lies in a bounded set of $L^2(\Sigma_0)\times L^2(\Sigma_0)$. We argue by contradiction. We assume that there exists four sequences $(a_n)_{n\in\mathbb N}$, $(b_n)_{n\in\mathbb N}$, $(\gamma_{1,n})_{n\in\mathbb N}$ and  $(\gamma_{2,n})_{n\in\mathbb N}$ such that
\begin{itemize}
\item $a_n\cvfstar a\text{ in }L^\infty(Q_0)$;
\item $b_n\cvfstar b\text{ in }L^\infty(Q_0)$;
\item $\gamma_{1,n}\rightarrow \gamma_1\text{ in } L^2(\Sigma_0)$;
\item  $\gamma_{2,n}\rightarrow \gamma_2\text{ in } L^2(\Sigma_0)$;
\item $\Vert \delta_{1,n}^\varepsilon\Vert _{L^2(\Sigma_0)}^2+\Vert \delta_{2,n}^\varepsilon\Vert _{L^2(\Sigma_0)}^2\rightarrow \infty,$ where the pair $(\delta_{1,n}^\varepsilon,\delta_{2,n}^\varepsilon)$ denotes the unique minimizer of $J_\varepsilon$ (for the potentials $a_n$ and $b_n$ instead of $a$ and $b$, and boundary conditions $\gamma_{1,n}$ and $\gamma_{2,n}$ instead of $\delta_1$ and $\delta_2$, in system \eqref{adj1}-\eqref{adj2}).
\end{itemize}
 Let us show that this implies 
\begin{equation}\label{Jn}J_\varepsilon (\delta_{1,n}^\varepsilon,\delta_{2,n}^\varepsilon)\rightarrow \infty, \quad \text{ as  }n\rightarrow\infty,\end{equation}
which will lead to a contradiction, since according to the point (i), we should have $J(\delta_{1,n},\delta_{2,n})\leqslant 0$ for every $n\in\mathbb N$.
We introduce 
$$
\tilde\delta_{1,n}^\varepsilon:=\frac{\delta_{1,n}^\varepsilon}{\displaystyle\sqrt{
\Vert \delta_{1,n}^\varepsilon\Vert _{L^2(\Sigma_0)}^2+\Vert \delta_{2,n}^\varepsilon\Vert_{L^2(\Sigma_0)}^2}}
$$ 
and 
$$
\tilde\delta_{2,n}^\varepsilon:=\frac{\delta_{2,n}^\varepsilon}{\displaystyle\sqrt{
\Vert \delta_{1,n}^\varepsilon \Vert_{L^2(\Sigma_0)}^2+\Vert \delta_{2,n}^\varepsilon\Vert_{L^2(\Sigma_0)} ^2}}.$$
One clearly has that $(\tilde \delta_{1,n}^\varepsilon)$ and $(\tilde\delta_{2,n}^\varepsilon)$ are bounded in $L^2(\Sigma_0)$, so that we may assume without loss of generality that they converge weakly respectively to $\tilde \delta_{1,opt}^\varepsilon\in L^2(\Sigma_0)$ and $\tilde \delta_{2,opt}^\varepsilon\in L^2(\Sigma_0)$. Hence, a compact embedding argument allows to prove that the corresponding solution $(\tilde \varphi_n^\varepsilon,\tilde\psi_n^\varepsilon)$ of \eqref{adj1}-\eqref{adj2} converges strongly in $(L^2(Q_0))^2$ respectively  to some $(\tilde \varphi^\varepsilon,\tilde\psi^\varepsilon)$ which is still a solution of \eqref{adj1}-\eqref{adj2} with, as boundary term, $(\tilde \delta_{1,opt}^\varepsilon,\tilde\delta_{2,opt}^\varepsilon)$.

According to the point (i), we have
\begin{eqnarray*}
J_\varepsilon (\delta_{1,n}^\varepsilon,\delta_{2,n}^\varepsilon)&=& \left(\Vert \delta_{1,n}^\varepsilon\Vert^2+\Vert \delta_{2,n}^\varepsilon\Vert^2\right)\frac{1}{2}\int_{(0,T)\times\omega} (\tilde\psi^\varepsilon_n)^2+\varepsilon\displaystyle\sqrt{\Vert\delta_{1,n}^\varepsilon\Vert^2_{L^2(\Sigma_0)}+\Vert\delta_{2,n}^\varepsilon\Vert^2_{L^2(\Sigma_0)}}\\
&& +\displaystyle\sqrt{\left(\Vert\delta_{1,n}^\varepsilon\Vert^2_{L^2(\Sigma_0)}+\Vert\delta_{2,n}^\varepsilon\Vert^2_{L^2(\Sigma_0)}\right )}\left(-\int_{\Sigma_0}\gamma_{1,n}\tilde\delta_{1,n}^\varepsilon+\int_{\Sigma_0}\gamma_{2,n}\tilde\delta_{2,n}^\varepsilon\right )\\&
\leqslant 0.
\end{eqnarray*}
Hence, dividing each side of this inequality by $\Vert\delta_{1,n}^\varepsilon\Vert^2_{L^2(\Sigma_0)}+\Vert\delta_{2,n}^\varepsilon\Vert^2_{L^2(\Sigma_0)}$, we infer that 
$$
\int_{(0,T)\times\omega} |\tilde \psi_n^\varepsilon|^2\rightarrow 0\quad \text{ as }n\rightarrow\infty, \quad \text{ i.e. } \,\mbox{ }\,\tilde\psi^\varepsilon=0 \text{ on }(0,T)\times\omega.
$$ 
By using the same unique continuation argument as in the proof of Proposition \ref{mainpr}, we infer that $\tilde\delta_{1,opt}^\varepsilon=\tilde\delta_{2,opt}^\varepsilon=0$.
Now, since there holds
$$
\frac{J_\varepsilon(\delta_{1,n}^\varepsilon,\delta_{2,n}^\varepsilon)}{\sqrt{\Vert\delta_{1,n}^\varepsilon\Vert^2_{L^2(\Sigma_0)}+\Vert\delta_{2,n}^\varepsilon\Vert^2_{L^2(\Sigma_0)}}}\geqslant\varepsilon\displaystyle-\int_{\Sigma_0}\gamma_{1,n}\tilde\delta_{1,n}^\varepsilon+\int_{\Sigma_0}\gamma_{2,n}\tilde \delta_{2,n}^\varepsilon,
$$
it follows that \eqref{Jn} holds true since the right-hand side converges to $\varepsilon>0$ as $n\rightarrow\infty$.

As a consequence, the pair $(\delta_{1,opt}^\varepsilon,\delta_{2,opt}^\varepsilon)$ lies in a bounded set of $L^2(\Sigma_0)\times L^2(\Sigma_0)$, from which we deduce with the help of a standard compact embedding argument that $k_{opt}^\varepsilon=\psi_{opt}^\varepsilon\chi_\omega$ lies in a compact set of $L^2(Q_0)$.
\subsection{The semi-linear case (Proof of Theorem \ref{th1})}
Let us go back to the Proof of Theorem \ref{th1}. The proof below is very similar to \cite[Proof of Theorem 1]{2}, so that we will skip some details and recall only the main lines.

Let us introduce
\begin{equation*}
F(s):=\frac{f(s)-f(0)}{s}.\end{equation*}
Then $F$ is continuous and bounded on $\mathbb R$ since $f\in C^1(\R)\cap \operatorname{Lip}(\R)$, where $\operatorname{Lip}(\R)$ denotes the set of Lipschitz functions on $\R$.
Let $z\in L^2(Q_0)$. We consider the linear system

\begin{equation}
\left\{ 
\begin{array}{rcll}
\dfrac{\partial u}{\partial t}-\Delta u +F(z)u&=&k\chi _{\omega
} & \text{ in }Q_{0} ,\\ 
u&=& 0 & \text{ on }\Sigma_{0},\\ 
u\left( 0,.\right) &=& 0 & \text{ in }\Omega_{0},
\end{array}%
\right.\label{pf1}
\end{equation}
\begin{equation}
\left\{ 
\begin{array}{rcll}
-\dfrac{\partial v}{\partial t}-\Delta v +f'(z)v &=& u\chi_{\Theta}& \text{ in }Q_{0} ,\\ 
v&=& 0 & \text{ on }\Sigma_{0},\\ 
v\left( T,.\right) &=& 0 & \text{ in }\Omega_{0},
\end{array}
\right.\label{pf2}
\end{equation}
\begin{equation}
\left\{ 
\begin{array}{rcll}
\dfrac{\partial u_0}{\partial t}-\Delta v_0 +F(z)u_0 &=&-f(0)+\xi & \text{ in }Q_{0} ,\\ 
u_0&=& 0 & \text{ on }\Sigma_{0},\\ 
u_0\left( 0,.\right) &=& 0 & \text{ in }\Omega_{0},
\end{array}%
\right.\label{pf3}
\end{equation}
\begin{equation}
\left\{ 
\begin{array}{rcll}
-\dfrac{\partial v_0}{\partial t}-\Delta v_0 +f'(z)v_0 &=&u_0\chi_{\Theta} & \text{ in }Q_{0},\\ 
v_0&=& 0 & \text{ on }\Sigma_{0},\\ 
v_0\left( T,.\right) &=&0 & \text{ in }\Omega_{0},
\end{array}%
\right.\label{pf4}
\end{equation}

Let $\varepsilon>0$ and $k_{opt}^\varepsilon$ be the optimal control defined in Proposition \ref{mainpr} with $a=F(z)$, $b=f'(z)$, $\gamma_1=\partial_n u_0$ and $\gamma_2=\partial_n v_0$, and we still call by $(u,v)$ the corresponding solution to \eqref{pf1}-\eqref{pf2} for the sake of simplicity.
We now introduce
\begin{equation}\label{pfy}y(t,x)=u(t,x)+u_0(t,x) \end{equation}
and 
\begin{equation*}
q(t,x)=v(t,x)+v_0(t,x),\end{equation*}
as well as the nonlinear operator 
\begin{equation}\label{lambda}
\Lambda:z\in L^2(Q_0)\longmapsto y\in L^2(Q_0).
\end{equation}
It is standard that any fixed point of $\Lambda$ will provide a solution $(z,q)$ of \eqref{cou1}-\eqref{cou2} satisfying moreover 
	$$
	\displaystyle\sqrt{\int_{\Sigma_0} |\partial_{n} z|^2+\int_{\Sigma_0}|\partial_{n} q |^2}\leqslant \varepsilon,
	$$
so that \eqref{afin} holds true, as a consequence of the Cauchy-Schwarz and Young inequalities.

We remark that $F(L^2(Q_0))$ and $f'(L^2(Q_0))$ are bounded sets of $L^\infty(Q_0)$ and that the solution $(u,v)$ of \eqref{pf1}-\eqref{pf2} depends continuously on their data, so that $\Lambda$ is continuous.  

The expected result will be derived by applying the Schauder fixed point theorem. 

According to the considerations above, it remains to prove that the range of $\Lambda$ is relatively compact in $L^2(Q_0)$ (here and in the sequel, this wording must be understood for the strong topology of $L^2$). We follow \cite[Pages 677-678]{2}. By using standard arguments from variational analysis, when $z$ spans $L^2(Q_0)$, $f'(z)$ lies in a bounded set of $L^\infty(Q_0)$, so that the solution $u_0$ of \eqref{pf3} lies in a bounded set in $L^2((0,T),H^2(\Omega_0))$. Hence,  we infer that $\partial_n u_0$ lies in a relatively compact set in $L^2(\Sigma_0)$. The same argument applies for the solution $v_0$ of \eqref{pf4} since $F(z)$ lies in a bounded set of  $L^\infty(Q_0)$ and $u_0\chi_{\Theta}$ in a bounded (and even relatively compact) set in $L^2(Q_0)$. Hence, the pair $(\partial_n u_0,\partial_n v_0)$ lies in a relatively compact set in $L^2(\Sigma_0)\times L^2(\Sigma_0)$.

Applying point (iii) of Proposition \ref{mainpr}, we infer that the control $k_{opt}^\varepsilon$ lies in a bounded set of $L^2(\Sigma_0)$, so that the corresponding solution $u$ of \eqref{pf1} lies in a compact set of $L^2(Q_0)$. Since the solution $u_0$ of \eqref{pf3} lies in a bounded set in $L^2((0,T),H^2(\Omega_0))$, $u_0$ also lies in a relatively compact set in $L^2(Q_0)$. Hence, $y$ defined by \eqref{pfy} also lies in a relatively compact set of $L^2(Q_0)$, so that there exists a fixed point to the operator $\Lambda$ defined in \eqref{lambda} and the proof of Theorem \ref{th1} is complete.
\cqfd

\section{Insensitizing control for a one or two-parameter family of deformations in the linear case}\label{S2}

The exact insensitizing problem \eqref{efin}, although interesting, appears intricate and we did not manage to determine sufficient conditions on the parameters of the problem \eqref{nlh} allowing to ensure it, even in the linear case. In this section, we consider a simplified problem, looking at insensitizing the functional $\Phi$ for several subfamilies of deformations.
\subsection{Setting of the problem and comments on the main result}\label{sec:settingExactInsens}
To make the framework precise, consider some $M\in\mathbb N^*$ and a family of deformation 
$$
\mathcal{E}=\operatorname{span}\{\mathbf{V}_1,\dots,\mathbf{V}_M\},
$$
where $\mathbf{V}_i\in \mathcal{V}^{3,\infty}$ for all $i=1,\dots,M$.

The condition \eqref{eiFinite} can be recast in a much simpler form. Indeed, following the computations done in the proof of Proposition \ref{propcou}, this condition writes
\begin{equation*}
	\int_{\partial\Omega_{0}}(\mathbf{V}\cdot \mathbf{n})\left(\int_{0}^{T}\partial_{n} y_{\Omega_0}\partial_{n} q_{\Omega_0}\, dt \right)\, d\sigma=0,\, \quad \forall \mathbf{V}\in\mathcal{E},
\end{equation*}
where $(y_{\Omega_0},q_{\Omega_0})$ denotes the unique solution of the coupled system  \eqref{cou1}-\eqref{cou2}, or equivalently
\begin{equation}\label{eq:metz1304}
\mathcal{U}(h)=0\quad \text{where}\quad	[\mathcal{U}(h)]_i=\int_{\partial\Omega_{0}}(\mathbf{V}_i\cdot \mathbf{n})\left(\int_{0}^{T}\partial_{n} y_{\Omega_0}\partial_{n} q_{\Omega_0}\, dt \right)\, d\sigma,\, \quad i=1,\dots,M.
\end{equation}

The next section is devoted to proving Theorem \ref{th2}. Let us provide hereafter some examples of applications for particular choices of families $\mathcal{E}$.

\paragraph{Insensitizing with respect to translations/rotations in a plane.} 
Assume that $N\geqslant 2$. Insensitizing $\Phi$ with respect to all translations comes to consider the family $\mathcal{E}=\operatorname{span}\{\mathbf{\varepsilon}_1,\mathbf{\varepsilon}_2,\ldots,\mathbf{\varepsilon}_N\}$, with the $N$ (constant) vector fields 
$$\mathbf{\varepsilon}_1=\begin{pmatrix}1\\0\\0 \\ \vdots\\ 0\end{pmatrix},\ \mathbf{\varepsilon}_2=\begin{pmatrix}0\\ 1\\0\\ \vdots\\ 0 \end{pmatrix}, \ \ldots, \mathbf{\varepsilon}_N=\begin{pmatrix}0\\ 0 \\ \vdots\\0\\1\end{pmatrix}.$$ 
If we consider a plane $\mathcal P$ in $\mathbb R^N$, which is generated by two vectors $\mathbf{E}_1$ and $\mathbf{E}_2$, we can also restrict to all translations that are in the direction of this plane, leading to solve
\begin{equation*}
\int_{\partial\Omega_{0}}(\mathbf{E}_i\cdot \mathbf{n})\left(\int_{0}^{T}\partial_{n} y_{\Omega_0}\partial_{n} q_{\Omega_0}\, dt \right)\, d\sigma=0,\, \quad i=1,2.
\end{equation*}

Concerning now the rotations in a plane, assume for the sake of simplicity that we consider the set of rotations in the plane $\mathbf P=\operatorname{span} \{\mathbf{\varepsilon}_1,\mathbf{\varepsilon}_2\}$.
Notice first that any rotation $R_{x_0,y_0,\theta}$ in $\mathbf P$ (assumed to be extended by the identity on $\mathbf P^\perp$), parametrized by $(x_0,y_0,\theta)\in \R^2\times \mathbb{S}^1$ where $(x_0,y_0)$ denote the coordinates of the rotation center and $\theta$ its angle, is given by
$$
R_{x_0,y_0,\theta}\begin{pmatrix}x\\y \end{pmatrix}=\begin{pmatrix}
x'=x_0+\cos\theta (x-x_0)-\sin \theta (y-y_0)\\
y'=y_0+\sin \theta(x-x_0)+\cos \theta (y-y_0)
\end{pmatrix}.
$$
Therefore, insensitizing with respect to all rotations of center $(x_0,y_0)$ in the plane $\mathbf P$ leads to consider the vector field $\mathbf{V}_{x_0,y_0}$ given by
$$
\mathbf{V}_{x_0,y_0}\begin{pmatrix}x\\y \end{pmatrix}=\lim_{\theta\to 0}\frac{1}{\theta}\left(R_{x_0,y_0,\theta}\begin{pmatrix}x\\y \end{pmatrix}-\begin{pmatrix}x\\y \end{pmatrix}\right)=\begin{pmatrix}-(y-y_0)\\x-x_0 \end{pmatrix}.
$$ 
As a consequence, insensitizing with respect to all rotations in the plane leads to consider the family $\mathcal{E}=\operatorname{span}\{\mathbf{V}_1,\mathbf{V}_2,\mathbf{V}_3\}$ with
$$
\mathbf{V}_1=\begin{pmatrix}1\\0\\0\\ \vdots\\ 0\end{pmatrix}, \quad \mathbf{V}_2=\begin{pmatrix}0\\1\\0\\\vdots\\ 0\end{pmatrix}\quad \text{and}\quad
\mathbf{V}_3=\begin{pmatrix}-y\\x\\0\\ \vdots\\ 0\end{pmatrix}.
$$
Our results does not cover the case of families of deformation of dimension 3. In the case where $\Omega_0$ is a cylinder of the form $\mathcal B_2(0,R)\times \prod_{i=3}^N(a_i,b_i)$ (where $\mathcal B_2(0,R)$ is the two-dimensional euclidean ball of radius $R>0$ and $a_i<b_i$), one can notice that $\mathbf{V}_3\cdot \mathbf{n}=0$ on $\Sigma_0$. Therefore, these considerations lead to the following byproduct of Theorem \ref{th2}.

\begin{corollary}
For any bounded connected domain $\Omega_0$ of class $C^2$, there exists a control $h_T$ insensitizing $\Phi$ at the first order with respect to all translations in a plane. Furthermore, if $\Omega_0$ is a cylinder of the form $\mathcal B_2(0,R)\times \prod_{i=3}^N(a_i,b_i)$ (where $\mathcal B_2(0,R)$ is the two-dimensional euclidean ball of radius $R>0$ and $a_i<b_i$), the control $h_T$ also insensistizes $\Phi$ at the first order with respect to all rotations of the plane $\operatorname{span} \{\mathbf{\varepsilon}_1,\mathbf{\varepsilon}_2\}$.
\end{corollary}


\subsection{Proof of Theorem \ref{th2} (case of one/two dimensional families of perturbations)}

The proof uses at the same time the density results stated in Proposition \ref{mainpr} as well as geometrical properties of second order curves.

\paragraph{Proof of the case $M=1$.} Let us consider that $\mathcal{E}=\operatorname{span}(\mathbf{V})$. Recall that, according to \eqref{eq:metz1304}, the problem comes to determine a control function $h$ such that
\begin{equation}\label{argenton}
\int_{\partial\Omega_{0}}(\mathbf{V}\cdot \mathbf{n})\left(\int_{0}^{T}\partial_{n} y_{\Omega_0}\partial_{n} q_{\Omega_0}\, dt \right)\, d\sigma=0.
\end{equation}
First of all, one remarks that if $\mathbf{V}\cdot \mathbf{n}=0$ on $\partial\Omega$ (which is possible since we only assumed that $\mathbf{V} $ is non-zero as a function defined on the whole space $\mathbb R^N$), then \eqref{argenton} is automatically verified. Hence, from now on we assume that $\mathbf{V}\cdot \mathbf{n}\not\equiv 0$ on $\partial\Omega$.

Let us recast this question. For that purpose, we introduce the pairs $(y_\xi,q_\xi)$ and $(y_h,q_h)$ as the solution of the linear systems
	\begin{equation}
\left\{ 
\begin{array}{rcll}
\dfrac{\partial y_\xi}{\partial t}-\Delta y_\xi &=&\xi & \text{ in }Q_{0} ,\\ 
y_\xi &=& 0 & \text{ on }\Sigma_{0},\\ 
y_\xi\left( 0,\cdot \right) &=& 0 & \text{ in }\Omega_{0},
\end{array}%
\right.\label{cou1b}
\end{equation}
\begin{equation}
\left\{ 
\begin{array}{rcll}
-\dfrac{\partial q_\xi}{\partial t}-\Delta q_\xi &=& y_\xi\chi_{\Theta}& \text{ in }Q_{0} ,\\ 
q_\xi&=& 0 & \text{ on }\Sigma_{0},\\ 
q_\xi\left( T,\cdot\right) &=& 0 & \text{ in }\Omega_{0},
\end{array}%
\right.\label{cou2b}
\end{equation}
and
	\begin{equation}
\left\{ 
\begin{array}{rcll}
\dfrac{\partial y_h}{\partial t}-\Delta y_h &=& h\chi _{\omega} & \text{ in }Q_{0} ,\\ 
y_h &=& 0 & \text{ on }\Sigma_{0},\\ 
y_h \left( 0,\cdot \right) &=& 0 & \text{ in }\Omega_{0},
\end{array}%
\right.\label{cou1bis}
\end{equation}
\begin{equation}
\left\{ 
\begin{array}{rcll}
-\dfrac{\partial q_h}{\partial t}-\Delta q_h  &=& y_h\chi_{\Theta}& \text{ in }Q_{0} ,\\ 
q_h &=& 0 & \text{ on }\Sigma_{0},\\ 
q_h\left( T,\cdot\right) &=& 0 & \text{ in }\Omega_{0},
\end{array}%
\right.\label{cou2bis}
\end{equation}
in such a way that $y_{\Omega_0}=y_\xi+y_h$ and $q_{\Omega_0}=q_\xi+q_h$.

Then, the quantity $\mathcal{U}(h)\in \R$ defined in \eqref{eq:metz1304} can be decomposed as
$$
\mathcal{U}(h)=Q(h)+L(h)+C,
$$
where
\begin{eqnarray*}
Q(h)&=& \int_{\Sigma_0} (\mathbf{V}\cdot \mathbf{n})\partial_{n} y_{h}\partial_{n} q_{h}\, dt d\sigma,\\
L(h)&=& \int_{\Sigma_0}(\mathbf{V}\cdot \mathbf{n}) (\partial_{n} y_{\xi}\partial_{n} q_{h}+\partial_{n} y_{h}\partial_{n} q_{\xi})\, dt d\sigma,\\
C&=& \int_{\Sigma_0} (\mathbf{V}\cdot \mathbf{n})\partial_{n} y_{\xi}\partial_{n} q_{\xi}\, dt d\sigma .
\end{eqnarray*}
Using this decomposition together with the facts that $Q$ is quadratic in $h$ and $L$ is linear in $h$, we claim that the problem comes to find a control function $h$ such that the algebraic equation $\mathcal{U}(\lambda h)=0$ with unknown the real number $\lambda$, also writing 
$$
\lambda^2 Q(h)+\lambda L(h)+C=0,
$$
has a real solution. Hence, the control function $\lambda h$, where $\lambda$ denotes a solution of the polynomial equation above, will solve \eqref{argenton}. We then infer that it is enough to choose $h$ such that the discriminant of this equation is positive, namely $L(h)^2-4Q(h)C>0$. Let us consider two functions $\gamma_1$ and $\gamma_2$ in $L^2(\Sigma_0)$ satisfying $\mathbf D(\gamma_1,\gamma_2)>0$, where
$$
\mathbf D(\gamma_1,\gamma_2)=\left(\int_{\Sigma_0} (\mathbf{V}\cdot \mathbf{n})(\partial_{n} y_{\xi}\gamma_2+\gamma_1\partial_{n} q_{\xi})\, dt d\sigma\right)^2-4C\left(\int_{\Sigma_0} (\mathbf{V}\cdot \mathbf{n})\gamma_1\gamma_2\, dt d\sigma\right).
$$
We claim that the existence of two such functions is straightforward as soon as $\mathbf{V}\cdot \mathbf{n}\neq 0$ almost everywhere on $\partial\Omega$.
Now, for a given positive number $\varepsilon$, according to Proposition \ref{mainpr}, there exists $h\in L^2((0,T)\times \omega)$ such that $\partial_{n} y_h=\gamma_1+\varepsilon_1$ and $\partial_{n} q_h=\gamma_2+\varepsilon_2$ with
$\left\Vert\varepsilon_1   \right \Vert_{L^2(\Sigma_0)}^2+\left\Vert\varepsilon_2   \right \Vert_{L^2(\Sigma_0)}^2\leqslant \varepsilon^2$. Moreover, using several times the triangle and the Cauchy-Schwarz inequalities leads to the estimate
$$
|L(h)^2-4Q(h)C-\mathbf D(\gamma_1,\gamma_2)|\leq c(\gamma_1,\gamma_2,\Vert \partial_{n} y_{\xi}\Vert_{L^2(\Sigma_0)},\Vert \partial_{n} q_{\xi}\Vert_{L^2(\Sigma_0)})\max_{i=1,2}\Vert\mathbf{V}\cdot \mathbf{n} \Vert _{L^\infty(\Sigma_0)}^i\varepsilon^i ,
$$
where $c(\gamma_1,\gamma_2,\Vert \partial_{n} y_{\xi}\Vert_{L^2(\Sigma_0)},\Vert \partial_{n} q_{\xi}\Vert_{L^2(\Sigma_0)})\in \R_+^*$ does not depend on $h$. Hence, we infer that it is possible to choose $\varepsilon>0$ small enough so that $L(h)^2-4Q(h)C>0$ and the desired conclusion follows.
\cqfd
\paragraph{Proof of the case $M=2$.} We generalize the approach used in the case $M=1$, by recasting the main issue to determining whenever two curves of degree 2 in the plane intersect. Let us first  assume that $\mathbf V_1.\mathbf{n}$ and $\mathbf V_2.\mathbf{n} $ are linearly dependent as functions defined on $\partial \Omega_0$ (this is possible since $\{\mathbf V_1.\mathbf{n},\mathbf V_2.\mathbf{n} \}$  is assumed to be linearly independent as functions defined on the whole space $\mathbb R^N$). Then, there exists some $\nu\in\mathbb R$ such that $\mathbf V_i.\mathbf{n}=\nu \mathbf V_j.\mathbf{n}$ for $(i,j)=(1,2)$ or $(i,j)=(2,1)$, so that we have returned to the previous situation (i.e. $M=1$), which has already been treated.
Hence, we assume from now on that $\mathbf V_1.\mathbf{n}$ and $\mathbf V_2.\mathbf{n}$ are linearly independent as functions defined on $\partial \Omega_0$.

Let us consider two functions $h_1$ and $h_2$ (that we will choose adequately in the sequel) and let us write $h=\lambda h_1+\mu h_2$ where $(\lambda,\mu)\in \R^2$ will also be chosen in the sequel. Hence, the vector $\mathcal{U}(h)\in \R^2$ can be decomposed as
$$
\mathcal{U}(h)=\begin{pmatrix}
\lambda^2A_1(h_1)+2\lambda\mu B_1(h_1,h_2)+\mu^2C_1(h_2)+\lambda D_1(h_1)+\mu E_1(h_2)+F_1\\
\lambda^2A_2(h_1)+2\lambda \mu B_2(h_1,h_2)+\mu^2C_2(h_2)+\lambda D_2(h_1)+\mu E_2(h_2)+F_2
\end{pmatrix},
$$
where, for $i=1,2$, one has
\begin{eqnarray*}
A_i(h_1)&=& \int_{\Sigma_0} (\mathbf{V}_i\cdot \mathbf{n})\partial_{n} y_{h_1}\partial_{n} q_{h_1}\, dt d\sigma,\\
B_i(h_1,h_2)&=& \frac{1}{2}\int_{\Sigma_0}(\mathbf{V}_i\cdot \mathbf{n}) (\partial_{n} y_{h_1}\partial_{n} q_{h_2}+\partial_{n} y_{h_2}\partial_{n} q_{h_1})\, dt d\sigma,\\
C_i(h_2) &=& \int_{\Sigma_0} (\mathbf{V}_i\cdot \mathbf{n})\partial_{n} y_{h_2}\partial_{n} q_{h_2}\, dt d\sigma,\\
D_i(h_1) &=& \int_{\Sigma_0} (\mathbf{V}_i\cdot \mathbf{n})(\partial_{n} y_{h_1}\partial_{n} q_{\xi}+\partial_{n} y_{\xi}\partial_{n} q_{h_1})\, dt d\sigma,\\
E_i(h_2) &=& \int_{\Sigma_0} (\mathbf{V}_i\cdot \mathbf{n})(\partial_{n} y_{h_2}\partial_{n} q_{\xi}+\partial_{n} y_{\xi}\partial_{n} q_{h_2})\, dt d\sigma,\\
F_i &=& \int_{\Sigma_0} (\mathbf{V}_i\cdot \mathbf{n})\partial_{n} y_{\xi}\partial_{n} q_{\xi}\, dt d\sigma,
\end{eqnarray*}
with $(y_{h_i},q_{h_i})$ the pair solving the coupled system \eqref{cou1bis}-\eqref{cou2bis} where $h$ has been replaced by $h_i$, and $(y_{\xi},q_{\xi})$, the pair solving the coupled system \eqref{cou1b}-\eqref{cou2b}.

First of all, let us exclude several trivial cases: If $F_i=0$ for $i=1$ and/or $i=2$, one chooses $h_i=0$ and use the previous result for $M=1$. Hence, we assume from now on that $F_i\not =0$ ($i=1,2$), meaning in particular that the function $(\mathbf{V}_i\cdot \mathbf{n})\partial_{n} y_{\xi}\partial_{n} q_{\xi}$ does not vanish identically on $\partial\Omega_0$ for $i=1,2$.

Our strategy in the case $M=2$ is to determine first a favorable choice of Neumann traces $(\partial_n y_{h_i},\partial_n q_{h_i})$ on $\partial\Omega$ for the solutions of System  \eqref{cou1bis}-\eqref{cou2bis} allowing to insensitize exactly $\Phi$ for the family $\mathcal{E}$. Hence, we will use a perturbative argument to prove that such a favorable choice is reachable with the help of two control functions $h_1$ and $h_2$. For this last step, we will strongly exploit Theorem \ref{th1} (and in particular the density result stated in Proposition \ref{mainpr}).

To be more precise, our reasoning can be split into two steps:
\begin{itemize}
\item[(i)] we assume that there exists two control functions $h_1$ and $h_2$ in $L^2(Q_0)$ such that $\partial_n q_{h_1}=\partial_n y_{h_2}=0$, $\partial_n y_{h_1}=f$ and $\partial_n q_{h_2}=g$ for some functions $f$ and $g$ in $L^2(\Sigma_0)$. We will then show that $f$ and $g$ can be chosen adequately to guarantee the existence of a pair $(\lambda,\mu)\in \R^2$ so that $\mathcal{U}(h)=0$;
\item[(ii)] according to Proposition \ref{mainpr}, we consider $\varepsilon>0$ and two control functions $h_{\varepsilon,i}$, $i=1,2$, such that 
\begin{equation}\label{cond1534lundi}
\left\Vert\partial_{n} q_{h_1}   \right \Vert_{L^2(\Sigma_0)}^2+\left\Vert\partial_{n} y_{h_2} \right \Vert_{L^2(\Sigma_0)}^2+\left\Vert\partial_{n} y_{h_1} -f  \right \Vert_{L^2(\Sigma_0)}^2+\left\Vert\partial_{n} q_{h_2}-g   \right \Vert_{L^2(\Sigma_0)}^2\leqslant \varepsilon^2.
\end{equation} 
We will prove the existence of a pair $(\lambda,\mu)\in \R^2$ such that $\mathcal{U}(h_{\varepsilon})=0$, with $h_{\varepsilon}=\lambda h_{\varepsilon,1}+\mu h_{\varepsilon,2}$, whenever $\varepsilon>0$ is chosen small enough.
\end{itemize}


\paragraph{Step (i). A favorable situation.} 
Assume that there exists two control functions $h_1$ and $h_2$ in $L^2(Q_0)$ such that $\partial_n q_{h_1}=\partial_n y_{h_2}=0$, $\partial_n y_{h_1}=f$ and $\partial_n q_{h_2}=g$ where $f$ and $g$ denote two elements in $L^2(\Sigma_0)$ that will be chosen in the sequel. In that case, one has
$$
\mathcal{U}(h)=\begin{pmatrix}
2\lambda\mu\widehat B_1+\lambda\widehat D_1+\mu\widehat E_1+F_1\\
2\lambda\mu\widehat B_2+\lambda\widehat D_2+\mu\widehat E_2+F_2
\end{pmatrix},
$$
where, for $i=1,2$, one has by definition
\begin{equation}
\label{leshat}
\widehat B_i=\frac{1}{2}\int_{\Sigma_0}(\mathbf{V}_i\cdot \mathbf{n}) fg\, dt d\sigma, \quad \widehat D_i=\int_{\Sigma_0} (\mathbf{V}_i\cdot \mathbf{n}) f\partial_{n} q_{\xi}\, dt d\sigma, \quad \widehat E_i=\int_{\Sigma_0} (\mathbf{V}_i\cdot \mathbf{n})\partial_{n} y_{\xi}g\, dt d\sigma.
\end{equation}
We introduce the following (possibly degenerated) hyperbolae
$$(\mathcal H_1)=\{(\lambda,\mu)|2\lambda\mu\widehat B_1+\lambda\widehat D_1+\mu\widehat E_1+F_1=0\}$$ 
and
$$(\mathcal H_2)=\{(\lambda,\mu)|2\lambda\mu\widehat B_2+\lambda\widehat D_2+\mu\widehat E_2+F_2=0\}.$$
Moreover, we introduce, for $i=1,2$,
$$(\mathcal H_i)^+=\{(\lambda,\mu)|2\lambda\mu\widehat B_i+\lambda\widehat D_i+\mu\widehat E_i+F_i>0\}$$ 
and
$$(\mathcal H_i)^-=\{(\lambda,\mu)|2\lambda\mu\widehat B_i+\lambda\widehat D_i+\mu\widehat E_i+F_i<0\}.$$

The following result is the main ingredient in order to understand how to choose  $f$ and $g$.
\begin{lemma}\label{lem:1622lundi}
Assume that for $i=1,2$, there holds
\begin{equation}\label{eq1448lundi}
\widehat B_i\neq 0,\quad \widehat E_i\widehat D_i\neq 2\widehat B_i F_i\quad \text{and}\quad \frac{\widehat D_2}{\widehat B_2}\neq \frac{\widehat D_1}{\widehat B_1}.
\end{equation}
Then, the two hyperbolae $(\mathcal H_1)$ and $(\mathcal H_2)$ are non-degenerate and they intersect 
in the $(\lambda,\mu)$-plane if and only if 
\begin{equation}\label{trainmetz2041}
\Delta(\widehat B_1,\widehat B_2,\widehat D_1,\widehat D_2,\widehat E_1,\widehat E_2)\leq 0,
\end{equation}
where
\begin{eqnarray*}
\lefteqn{\Delta(\widehat B_1,\widehat B_2,\widehat D_1,\widehat D_2,\widehat E_1,\widehat E_2)=\left(\frac{\widehat E_2}{2\widehat B_2}-\frac{\widehat E_1}{2\widehat B_1}\right)\left(\frac{\widehat D_2}{2\widehat B_2}-\frac{\widehat D_1}{2\widehat B_1}\right)}\phantom{aaaaa}\nonumber  \\
& & -\frac{1}{4}\left(\frac{F_1}{2\widehat B_1}-\frac{F_2}{2\widehat B_2}+\frac{\widehat E_1F_1}{4\widehat B_1^2}-\frac{\widehat E_2F_2}{4\widehat B_2^2}-\left(\frac{\widehat E_2}{2\widehat B_2}-\frac{\widehat E_1}{2\widehat B_1}\right)\left(\frac{\widehat D_2}{2\widehat B_2}-\frac{\widehat D_1}{2\widehat B_1}\right)\right)^2.
\end{eqnarray*}
Moreover, if
\begin{equation}\label{deltastr}
\Delta(\widehat B_1,\widehat B_2,\widehat D_1,\widehat D_2,\widehat E_1,\widehat E_2)<0,
\end{equation}
there are exactly two intersecting points and
\begin{equation}\label{labcn} (\mathcal H_1)\cap(\mathcal H_2)^+\not =\emptyset \mbox{ and }(\mathcal H_1)\cap(\mathcal H_2)^-\not =\emptyset,
\end{equation}
meaning that the two hyperbolae intersect non-tangentially.
\end{lemma}
\noindent\textbf{Proof of Lemma \ref{lem:1622lundi}.}
Noting that $\Omega_1(-\widehat E_1/(2\widehat B_1),-\widehat D_1/(2\widehat B_1))$ is the center of $(\mathcal{H}_1)$, we introduce the change of coordinates $U=\lambda+ \widehat E_1/(2\widehat B_1)$ and $V=\mu+\widehat D_1/(2\widehat B_1)$ and we will recast the equations of the hyperbolae in terms of the new coordinates $(U,V)$. This way, $(\mathcal{H}_1)$ becomes centered and its cartesian equation in the $(U,V)$-plane is
\begin{equation}\label{H1C}
UV=k_1\quad \text{with}\quad k_1=\frac{\widehat E_1\widehat D_1}{4\widehat B_1^2}-\frac{F_1}{2\widehat B_1}.
\end{equation}
Similarly, the equation of $(\mathcal{H}_2)$ in the $(U,V)$-plane is
\begin{equation}\label{H2C}
(U-u_2)(V-v_2)=k_2\quad \text{with}\quad u_2=\frac{\widehat E_2}{2\widehat B_2}-\frac{\widehat E_1}{2\widehat B_1}, \ v_2=\frac{\widehat D_2}{2\widehat B_2}-\frac{\widehat D_1}{2\widehat B_1}, \ k_2=\frac{\widehat E_2\widehat D_2}{4\widehat B_2^2}-\frac{F_2}{2\widehat B_2}.
\end{equation}
Since $k_1\not=0$ and $k_2\not =0$ by assumption, we infer that the two hyperbolas $(\mathcal H_1)$ and $(\mathcal H_2)$ are non-degenerate. Moreover, they intersect if and only if the system
$$
UV=k_1\quad \text{and}\quad (U-u_2)(V-v_2)=k_2
$$
has a solution. Plugging the relation $V=k_1/U$ into the relation $(U-u_2)(V-v_2)=k_2$ yields that the previous system has a solution if and only if the second order polynomial $v_2X^2+(k_2-k_1-u_2v_2)X+k_1u_2$ has a nonzero real root. Note that the assumption \eqref{eq1448lundi} yields in particular that $v_2\neq 0$. This is equivalent to claiming that the discriminant of this polynomial is nonnegative. It rewrites

$$u_2v_2k_1\leq \frac{1}{4}(k_2-k_1-u_2v_2)^2,$$
which is equivalent to condition \eqref{trainmetz2041}.
Assume now that \eqref{labcn} holds, which is equivalent to assume that   the polynomial $v_2X^2+(k_2-k_1-u_2v_2)X+k_1u_2$ has two distinct roots denoted by $r_1<r_2$.  It means that $(\mathcal H_1)$ and $(\mathcal H_2)$ have exactly two intersecting points. Consider some $V_1\in\mathbb R$ verifying $(V_1-r_1)(V_1-r_2)<0$. Then $(k_1/V_1,V_1)\in (\mathcal H_1)$ by \eqref{H1C}. Moreover, by \eqref{H2C} we have that $V_1\in (\mathcal H_2)^{\pm}$. Now, consider some $V_2\in\mathbb R$ verifying $(V_2-r_1)(V_2-r_2)>0$. Then $(k_1/V_2,V_2)\in (\mathcal H_1)$ by \eqref{H1C}. Moreover, by \eqref{H2C} we have that $V_2\in (\mathcal H_2)^{\mp}$.
The desired result follows.
\cqfd

We can now state the main result of this step.

\begin{lemma}\label{lemma:1557}
With the previous notations and under the assumptions of Theorem \ref{th2}, there exists $(f,g)\in (L^2(\Sigma_0))^2$ such that the two hyperbolae having for respective equations
$$
2\lambda\mu\widehat B_1+\lambda\widehat D_1+\mu\widehat E_1+F_1=0,\quad\text{and}\quad
2\lambda\mu\widehat B_2+\lambda\widehat D_2+\mu\widehat E_2+F_2=0,
$$
in the $(\lambda,\mu)$-plane, verifies \eqref{eq1448lundi} and \eqref{deltastr}.

\end{lemma}
\noindent \textbf{Proof of Lemma \ref{lemma:1557}.}
To simplify these conditions, we will use a homogeneity argument. Indeed, let $f$ and $g$ be fixed. Changing $f$ into $\eta f$ and $g$ into $\eta g$ with $\eta\in \R$ and making $\eta$ tend to $+\infty$ changes condition \eqref{eq1448lundi} into
\begin{equation}\label{cond:1502}
\widehat B_i\neq 0, \quad \widehat E_i\widehat D_i\neq 0,  \quad i=1,2,\quad \text{and}\quad \widehat D_2\widehat B_1-\widehat D_1\widehat B_2\neq 0,
\end{equation}
 and condition \eqref{deltastr} into
\begin{equation}\label{cond:1503}
\widehat F_2\widehat B_1-\widehat F_1\widehat B_2\neq 0.
\end{equation}

According to the expression given in \eqref{leshat}, and since  the product $\mathbf{V}_i\cdot \mathbf{n}\partial_{n} y_{\xi}\partial_{n} q_{\xi}$ ($i=1,2$) does not vanish identically by assumption, it is an easy task to construct two functions $f$ and $g$ such that \eqref{cond:1502} and \eqref{cond:1503} are verified. For instance, consider any $f\in L^2(\Sigma_0)$ such that, for $i=1,2$, 
\begin{eqnarray*}
\int_{\Sigma_0}(\mathbf{V}_i\cdot \mathbf{n}) f dt d\sigma\neq 0,\\
\int_{\Sigma_0} (\mathbf{V}_i\cdot \mathbf{n}) \partial_{n} q_{\xi}f\, dt d\sigma \not =0,\\
\int_{\Sigma_0} (\mathbf{V}_i\cdot \mathbf{n})\partial_{n} y_{\xi}f\, dt d\sigma\not = 0,\\
\int_{\Sigma_0}(\mathbf{V}_i\cdot \mathbf{n}) \partial_{n} q_{\xi}f^2 dt d\sigma\neq 0,\\
\end{eqnarray*}
and consider $g(\nu)=1+\nu f$. According to conditions \eqref{cond:1502} and \eqref{cond:1503}, it is enough to find some $\nu\in\mathbb R$ such that all the following conditions are verified, for $i=1,2$:

\begin{eqnarray*}
\int_{\Sigma_0}(\mathbf{V}_i\cdot \mathbf{n})  dt d\sigma+\nu\int_{\Sigma_0}(\mathbf{V}_i\cdot \mathbf{n}) f dt d\sigma\neq 0,\\\
\nu \left(\int_{\Sigma_0} (\mathbf{V}_1\cdot \mathbf{n}) \partial_{n} q_{\xi}f\, dt d\sigma \int_{\Sigma_0} (\mathbf{V}_2\cdot \mathbf{n}) \partial_{n} q_{\xi}f^2\, dt d\sigma \right . \\ \left .-\int_{\Sigma_0} (\mathbf{V}_2\cdot \mathbf{n}) \partial_{n} q_{\xi}f\, dt d\sigma \int_{\Sigma_0} (\mathbf{V}_1\cdot \mathbf{n}) \partial_{n} q_{\xi}f^2\, dt d\sigma \right )\not =0,
\\
\left(\int_{\Sigma_0} (\mathbf{V}_2\cdot \mathbf{n})\partial_{n} y_{\xi}\partial_{n} q_{\xi}\, dt d\sigma\right)\left(\int_{\Sigma_0} (\mathbf{V}_1\cdot \mathbf{n}) \partial_{n} q_{\xi}f+\nu\int_{\Sigma_0} (\mathbf{V}_1\cdot \mathbf{n}) \partial_{n} q_{\xi}f^2\, dt d\sigma \right )\\
- \left(\int_{\Sigma_0} (\mathbf{V}_1\cdot \mathbf{n})\partial_{n} y_{\xi}\partial_{n} q_{\xi}\, dt d\sigma \right)\left(\int_{\Sigma_0} (\mathbf{V}_2\cdot \mathbf{n}) \partial_{n} q_{\xi}f+\nu\int_{\Sigma_0}( \mathbf{V}_2\cdot \mathbf{n}) \partial_{n} q_{\xi}f^2\, dt d\sigma \right) 
\not =0.\\
\end{eqnarray*}
It is obvious that all real parameter $\mu$ apart from a finite number of values verifies the above relations and the result follows.
\cqfd
\paragraph{Step (ii). Use of the density result. }Let $\varepsilon>0$ and $(f,g)$ be chosen as in the statement of Lemma \ref{lemma:1557}.
According to Proposition \ref{mainpr}, we consider $\varepsilon>0$ and two control functions $h_{\varepsilon,i}$, $i=1,2$ such that the condition \eqref{cond1534lundi} holds true.
Let $h_{\varepsilon}=uh_{\varepsilon,1}+vh_{\varepsilon,2}$ for some $(\lambda,\mu)\in \R^2$. Then, one shows easily by using several times the Cauchy-Schwarz inequality that
$$
\mathcal{U}(h_{\varepsilon})=\begin{pmatrix}
u^2A_1(h_{\varepsilon,1})+2uvB_1(h_{\varepsilon,1},h_{\varepsilon,2})+v^2C_1(h_{\varepsilon,2})+uD_1(h_{\varepsilon,1})+vE_1(h_{\varepsilon,2})+F_1\\
u^2A_2(h_{\varepsilon,1})+2uvB_2(h_{\varepsilon,1},h_{\varepsilon,2})+v^2C_2(h_{\varepsilon,2})+uD_2(h_{\varepsilon,1})+vE_2(h_{\varepsilon,2})+F_2
\end{pmatrix},
$$
with $A_i(h_{\varepsilon,1})=\operatorname{o}(\varepsilon)$, $ B_i(h_{\varepsilon,1},h_{\varepsilon,2})=\widehat B_i+\operatorname{o}(\varepsilon)$, $C_i(h_{\varepsilon,2})=\operatorname{o}(\varepsilon)$, $D_i(h_{\varepsilon,1})=\widehat D_i+\operatorname{o}(\varepsilon)$, $E_i(h_{\varepsilon,2})=\widehat E_i+\operatorname{o}(\varepsilon)$ for $i=1,2$.  Moreover, we claim that, by construction, the coefficients $A_i(h_{\varepsilon,1})$, $ B_i(h_{\varepsilon,1},h_{\varepsilon,2})$, $C_i(h_{\varepsilon,2})$, $D_i(h_{\varepsilon,1})$, $E_i(h_{\varepsilon,2})$ for $i=1,2$ can be chosen as continuous functions of $\varepsilon$ in a neighborhood of 0.
 
 As a consequence, whenever $\varepsilon>0$ is chosen small enough, each line of the system $\mathcal{U}(h_{\varepsilon})=0$ defines a non-degenerated hyperbolic curve since $ A_i(h_{\varepsilon,1})C_i(h_{\varepsilon,2})-B_i(h_{\varepsilon,1},h_{\varepsilon,2})^2=\widehat A_i\widehat C_i-\widehat B_i^2+\operatorname{o}(\varepsilon)$ for $i=1,2$. Moreover, the eigenelements associated to the matrices of the quadratic form defining these hyperbolae converge to the eigenelements of the limit hyperbolae.
 
Finally, it also follows that for $\varepsilon>0$ small enough, the two hyperbolae $(\mathcal{H}_1^\varepsilon)$ and $(\mathcal{H}_2^\varepsilon)$  having for cartesian equation $[\mathcal{U}(h_{\varepsilon})]_1=0$ and $[\mathcal{U}(h_{\varepsilon})]_2=0$ in the $(\lambda,\mu)$-plane, converge uniformly to the respective hyperbolae $(\mathcal{H}_1)$ and $(\mathcal{H}_2)$ defined in Lemma \ref{lem:1622lundi} on each compact of $\R^2$. Since  $(\mathcal{H}_1)$ and $(\mathcal{H}_2)$ meet non-tangentially, it follows that $(\mathcal{H}_1^\varepsilon)$ and $(\mathcal{H}_2^\varepsilon)$ do not meet tangentially as soon as $\varepsilon$ is small enough.
\cqfd

\section{Conclusion}
\label{CON}
Considering a semi-linear heat equation with Dirichlet boundary conditions and globally Lipschitz nonlinearity, we investigated the issue of insensitizing a quadratic functional of the state with respect to domain variations. We have first proved an approximated insensitizing property and second an exact insensitizing property for some finite dimensional families of deformations, for this functional.

Some open issues and generalizations remain to be investigated. They are in order:

\begin{itemize}
\item \textbf{exact insensitizing property.} We did not conclude about the characterization of all domains $\Omega_0$ for which there exists a control $h$ insensitizing exactly the functional $\Phi$. We have no conjecture to formulate about this issue.
\item \textbf{about the generalization of Theorem \ref{th2}.} It is plausible  that the statement of Theorem \ref{th2} can be generalized (at least generically with respect to families of deformations in $\mathcal{V}^{3,\infty}$) to an arbitrary finite number $M>2$ of perturbations. We nevertheless did not manage to prove it. Indeed, our arguments for the low dimensional case rest upon the fact that, by using geometrical considerations, we were able to recast the issue to the one of determining a control $h$ satisfying an ``open'' condition, that is a condition of the kind $G(h)>0$ where $G$ is a functional enjoying some nice continuity properties. In higher dimensional cases, such a trick seems more intricate to apply.
\item \textbf{extension of our results to other functionals.} We foresee to investigate generalizations of our two main results to more general nonlinear functionals. In such a case, the insensitizing condition \eqref{efin} will involve the use of an adjoint term depending nonlinearly of the state $y_{\Omega_0}$, making the underlying mathematical analysis more intricate.
\item \textbf{extension of our results to other equations.} One may wonder if the same kind of results can be proved for larger classes of equations, notably of hyperbolic type (like the wave or Euler equations) or of dispersive type (like the Schrödinger equation). However, the situation is likely to be much more intricate because of some geometric conditions that may appear on the control domain $\omega$, that are not necessarily  stable under perturbations of $\Omega$, and are likely to be not separated from $\partial\Omega$.

\end{itemize}

\footnotesize{
}

\end{document}